\documentclass[12pt,a4paper]{article}
\usepackage[T2A]{fontenc}
\usepackage[utf8]{inputenc}
\usepackage[english,russian]{babel}
\usepackage[tbtags]{amsmath}
\usepackage{amsfonts,amssymb,mathrsfs,amscd}

\usepackage{srcltx}\usepackage{graphicx,color
}
\usepackage[all]{xy}
\usepackage{latexsym,amsfonts,amssymb,amsmath,longtable,
mathrsfs,theorem,cite
}
\usepackage{stmaryrd}
\usepackage[hyper]{amsbib}
\newcommand{\BS}{{
\mathcal B}{\mathcal S}}

\theoremstyle{plain}
\newtheorem{Lemma}{{\bfseries Лемма}}

\newtheorem{Theo}{{\bfseries Теорема}}

\theoremstyle{definition}
\newtheorem{proof}{Доказательство}
\begin{document}

\begin{center}
{\Large\bf On the sharp Baer--Suzuki theorem for $\pi$-radicals:\\ sporadic groups  
\footnote{Z.\,Wu was supported by the Natural Science Foundation of Jiangsu Province, China (No.BK20210442) and Jiangsu Shuangchuang (Mass Innovation and Entrepreneurship) Talent Program (No.JSSCBS20210841). D.\,O.\,Revin was supported by the Russian Science Foundation, grant~No~19-11-00039. }\\
}

~

~

{\bf Nanying Yang, Zhenfeng Wu, and Danila O. Revin
}
\end{center}

\def\abstractname{\bf Abstract}


\begin{abstract}
Let
$\pi$ be a proper subset of the set of all primes. Denote
by
$r$ the smallest prime which does not belong to
$\pi$ and set
$m = r$ if
$r = 2$ or
$3$ and
$m = r-1$ if
$r \geqslant  5$.
We study the following conjecture: a conjugacy class
$D$ of a finite group
$G$ is contained in the
$\pi$-radical
$\mathrm{O}_\pi(G)$ of
$G$ if and only if every
$m$ elements of
$D$ generate a
$\pi$-subgroup. We confirm this conjecture
for each group
$G$ whose nonabelian composition factors are isomorphic to sporadic or alternating groups.


\end{abstract}

\section*{Введение}

В работе рассматриваются только конечные группы, и мы под термином ``группа'' будем подразумевать конечную группу.
  Через всюду
  $\pi$
  обозначается некоторое множество простых чисел.
Группа называется
{\it
$\pi$-группой},
если все простые делители ее порядка
принадлежат~$\pi$.
Используются следующие стандартные обозначения. Для группы
$G$ через
$\mathrm{O}_\pi(G)$ обозначается
{\it
$\pi$-радикал},
т.~е. наибольшая нормальная
$\pi$-подгруппа
группы~$G$. Если
$M$~--- подмножество
группы~$G$, то через
$\langle M\rangle$
обозначена подгруппа, порожденная
множеством~$M$.

Теорема Бэра~---~Сузуки \cite{Baer,Suz,AlpLy} утверждает

\medskip\noindent
{\bf Теорема Бэра~---~Сузуки.}  {\it Пусть $p$ --- некоторое
простое число, $G$ ---  конечная группа и $x\in G$. Тогда
$x\in \mathrm{O}_p(G)$ если и только если  $\langle x,x^g \rangle$ является $p$-группой для
любого $g\in G$. }
\medskip

 Теорема Бэра~---~Су\-зу\-ки имеет несколько эквивалентных формулировок.
Обобщения и аналоги теоремы Бэра~---~Сузуки исследовались многими авторами, см., в~частности,~\cite{AlpLy,MOS,Mamont,Mamont1,Soz,Gu,FGG,GGKP,GGKP1,GGKP2,GuestLevy, Palchik, Tyut,Tyut1,BS_odd,BS_Dpi,YRV,YWRV}. Например, Н.~Гордеевым, Ф.~Грюневальдом, Б.~Кунявским и Е.~Плоткиным \cite{GGKP2} и независимо П.~Флавеллом, С.~Гэстом и Р.~Гуральником \cite{FGG} доказано, что если любые четыре элемента из данного класса сопряженности в конечной группе порождают разрешимую группу, то  весь класс содержится в разрешимом радикале группы.

В работе \cite[теорема~1]{YRV}  было доказано, что для любого множества $\pi$ простых чисел существует число $m=m(\pi)$ такое, что в произвольной группе $G$ выполнено равенство
$$
\mathrm{O}_\pi(G)=\{x\in G\mid \langle x^{g_1}, \dots, x^{g_m}\rangle~\text{--- } \pi\text{-группа для любых } g_1,\dots g_m\in G\}.
$$ Наименьшее такое число $m$ согласно \cite{GGKP1} называется {\it шириной Бэра~---~Сузуки} класса $\pi$-групп и обозначается $\mathrm{BS}(\pi)$. Кроме того, доказано \cite[теорема~2]{YRV}, что если $\pi$~--- непустое собственное подмножество множества всех простых чисел и $r$~--- наименьшее простое число, не  лежащее в~$\pi$, то
$$
r-1\leq \mathrm{BS}(\pi)\leq \max\{11, 2(r-2)\}.
$$
Там же было высказано предположение, что нижняя оценка $r-1$ в большинстве случаев совпадает с~$\mathrm{BS}(\pi)$. Более точно:

\medskip\noindent
{\bf Гипотеза 1} \cite{YRV}.  { Пусть $\pi$ --- собственное подмножество множества всех простых чисел, содержащее как минимум два простых числа, и $r$~--- наименьшее простое число, не входящее в~$\pi$. Тогда $$\mathrm{BS}(\pi)\leq\left\{\begin{array}{rl}
                               r, & \text{ если } r\in\{2,3\}, \\
                               r-1, &   \text{ если } r\geq 5.
                             \end{array}\right.
$$}

В.~Н.~Тютянов \cite{Tyut1} подтвердил гипотезу~1 при $r=2$.  Известны также примеры \cite[пример~2]{BS_odd} множеств $\pi$ с $r=3$ таких, что $\mathrm{BS}(\pi)>2$. Результаты \cite{BS_odd} сводят гипотезу к изучению т.~н. почти простых групп. Чтобы сформулировать утверждение о почти простых группах, которое гарантировало бы справедливость гипотезы 1, напомним обозначения, введенные в~\cite{GS} и~\cite{YRV}.

Пусть $L$~--- неабелева простая группа, $r$~--- простой делитель ее порядка, и ${x\in\mathrm{Aut}(L)}$~--- ее нетождественный автоморфизм. Мы отождествляем $L$ c подгруппой $\mathrm{Inn}(L)$ в~$\mathrm{Aut}(L)$. Cогласно~\cite{GS} через $\alpha(x,L)$  обозначено наименьшее количество $L$-сопряженных с $x$ элементов, которые порождают $\langle L,x\rangle$. По аналогии c $\alpha(x,L)$ в работе~\cite{YRV}  через $\beta_r(x,L)$  обозначено наименьшее количество $L$-сопряженных с $x$ элементов, которые порождают подгруппу группы $\langle L,x\rangle$, порядок которой делится\footnote{Если $r$ не делит $|L|$, то величина $\beta_r(x,L)$ не определена.} на~$r$. Ясно, что $\beta_r(x,L)\leqslant \alpha(x,L)$.

Учитывая имеющуюся редукцию к почти простым группам и результаты для $r=2$, см. \cite{Tyut1,BS_odd}, справедливость гипотезы 1 была бы доказана, если бы удалось установить справедливость следующего утверждения:

\medskip\noindent
{\bf Гипотеза 2} \cite{YRV}.
Пусть $L$~--- неабелева конечная простая группа, $x\in \mathrm{Aut}(L)$~--- элемент простого порядка. Пусть также $r$~--- нечетное простое число, и простой делитель  $s$ порядка группы $L$  таков, что ${s=r}$, если $r$ делит $|L|$, и $s> r$ в противном случае. Тогда
  $$
  \beta_{s}(x,L)\leq
  \left\{
  \begin{array}{rl}
  3,&\text{если }r=3,\\
  r-1,&\text{если }r>3.
  \end{array}
  \right.
  $$

Для случая, когда $L$~--- знакопеременная группа, гипотеза~2 верна \cite[предложение~2]{YRV}. В данной статье мы подтвердим гипотезу~2 для простых спорадических групп. Таким образом, справедлива

\begin{Theo}\label{main}
Пусть $L$~--- одна из $26$ спорадических групп, $x\in \mathrm{Aut}(L)$~--- элемент простого порядка. Пусть также $r$~--- нечетное простое число, и простой делитель  $s$ порядка группы $L$ таков, что $s=r$, если $r$ делит $|L|$, и $s> r$ в противном случае.  Тогда
  $$
  \beta_{s}(x,L)\leq
  \left\{
  \begin{array}{rl}
  3,&\text{если }r=3,\\
  r-1,&\text{если }r>3.
  \end{array}
  \right.
  $$
\end{Theo}

С помощью этой теоремы, используя результаты работ \cite{BS_odd,YRV}, получаем следующее утверждение.

\begin{Theo}\label{Cor}
Пусть $\pi$~--- некоторое множество простых чисел и $r$~--- наименьшее простое число, не лежащее в $\pi$. Положим
$$m=\left\{\begin{array}{rl}
                               r, & \text{ если } r\in\{2,3\}, \\
                               r-1, &   \text{ если } r\geq 5.
                             \end{array}\right.
$$
Тогда
$$
\mathrm{O}_\pi(G)=\{x\in G\mid \langle x^{g_1}, \dots, x^{g_m}\rangle~\text{--- } \pi\text{-группа для любых } g_1,\dots g_m\in G\}
$$
для любой конечной группы $G$, всякий неабелев композиционный фактор которой изоморфен спорадической или знакопеременной группе.
 \end{Theo}

\section{Предварительные результаты}

Следуя \cite{BS_odd}, используем следующее обозначение.  Пусть $\pi$  --- некоторое множество простых
чисел, $m$~--- неотрицательное целое число.   Для конечной группы $G$ будем писать $G\in{{\mathcal B}{\mathcal S}}_{\pi}^{m}$, если группа $G$ обладает следующим свойством: класс сопряженности $D$ группы $G$ тогда и только тогда содержится в~$\mathrm{O}_\pi(G)$, когда любые  $m$ элементов из $D$ порождают $\pi$-группу.

Для подмножества $\pi$ множества $\mathbb{P}$ всех простых чисел полагаем $\pi'=\mathbb{P}\setminus\pi$.

\begin{Lemma}\label{red} {\rm \cite[лемма 7]{BS_odd}} {\it Пусть $\mathcal{X}$~--- класс конечных групп, замкнутый относительно взятия нормальных подгрупп, гомоморфных образов и расширений и содержащий все $\pi$-группы. Допустим,  $\mathcal{X}\nsubseteq\BS_{\pi}^{m}$ для некоторого натурального $m\geq 2$, и группа ${G\in\BS_{\pi}^{m}\setminus\mathcal{X}}$ выбрана так, что ее порядок является наименьшим. Тогда группа $G$ содержит подгруппу $L$ и элемент $x$ такие, что
\begin{itemize}
\item[$(1)$] $L\trianglelefteqslant G$;
\item[$(2)$] $L$ является неабелевой простой группой;
\item[$(3)$] $L$ не является $\pi$- или $\pi'$-группой;
\item[$(4)$] $C_G(L)=1$;
\item[$(5)$] любые $m$ сопряженных с $x$ элементов порождают $\pi$-группу;
\item[$(6)$] $x$ имеет простой порядок, принадлежащий $\pi$;
 \item[$(7)$] $G=\langle x, L\rangle$.
\end{itemize}}
\end{Lemma}

\begin{Lemma} \label{2notinpi} {\rm \cite[теорема 1]{BS_odd}} {\it Если $2\not\in\pi$, то $\BS_\pi^2$ совпадает с классом всех групп.
}
\end{Lemma}

Следующая лемма очевидна.

\begin{Lemma}\label{estim}
Пусть $L$~--- неабелева конечная простая группа, $x,y\in \mathrm{Aut}(L)\setminus\{1\}$. Предположим, что $x\in \langle y^{g_1},\dots,y^{g_k}\rangle$ для некоторых $g_1,\dots,g_k\in L$. Тогда
 $$
 \beta_r(x,L)\leq k\cdot\beta_r(y,L)
 $$
 для любого простого делителя $r$ порядка группы $L$.
\end{Lemma}

\section{Доказательство теоремы~\ref{main}}

Для спорадических групп мы используем обозначения из~\cite{atlas}. Мы докажем теорему~\ref{main} в двух леммах~\ref{Inner} и~\ref{Outer}, отвечающих случаям, когда $x\in{\mathrm{Inn}(L)=L}$  и $x\in\mathrm{Aut}(L)\setminus\mathrm{Inn}(L)$. Доказательство основывается на использовании теории характеров и следующего результата, доказанного Л.~Ди Мартино, М.~А.~Пеллегрини и А.~Е.~Залесским.

\begin{Lemma}\label{DMPZ} {\rm\cite[теорема~3.1]{DMPZ}}
Пусть $L$~--- конечная спорадическая группа и $x$~--- неединичный элемент из $L$. Тогда имеют место следующие утверждения.
\begin{itemize}
  \item[$(1)$]  Пусть $L\ne M$ и $x\in L$~--- не инволюция. Тогда $\alpha(x,L)= 2$ за исключением следующих случаев:

\begin{itemize}
  \item[\rm{(a)}] $(L,x)\in\{ (J_2, 3a), (HS, 4a), (McL, 3a), (Ly, 3a), (Co_1, 3a),(Fi_{22}, 3a), (Fi_{23}, 3a),$ $ (Fi_{23}, 3b), (Fi_{24}', 3a), (Fi_{24}',3b)\}$.  В этих случаях $\alpha(x,L)= 3$;

\item[\rm{(b)}] $(L,x) = (Fi_{22}, 3b)$. В этом случае $2 \leqslant \alpha(x,L)\leqslant 3$;

\item[\rm{(c)}] $(L,x) = (Suz, 3a)$. В этом случае $3 \leqslant \alpha(x,L)\leqslant 4$.
\end{itemize}

\item[$(2)$]  Пусть $L\ne M$ и $x\in L$~---  инволюция. Тогда $\alpha(x,L)= 3$ за исключением следующих случаев:
\begin{itemize}
\item[\rm{(a)}] $(L,x) \in\{ (J_2, 2a), (Co_2, 2a), (B, 2a)\}$. В этих случаях  $\alpha(x,L)= 4$;

\item[\rm{(b)}] $(L,x) \in\{ (Fi_{22}, 2a), (Fi_{23}, 2a)\}$. В этих случаях $5 \leqslant \alpha(x,L)\leqslant 6$.
\end{itemize}

\item[$(3)$]  Если $L = M$ и  $x\in L$~--- не инволюция, то  $2 \leqslant \alpha(x,L)\leqslant 3$;

\item[$(4)$]  Если $L = M$ и  $x\in L$~--- инволюция, то  $3 \leqslant \alpha(x,L)\leqslant 4$.
\end{itemize}
\end{Lemma}

\begin{Lemma}\label{Inner}
Пусть $L$~--- одна из $26$ спорадических групп, $r\in\pi(L)$~--- нечетное число и $s\notin\pi(L)$~--- простое число. Тогда для любого элемента $x\in L$ простого порядка выполнены неравенства $$\beta_r(x,L)\leqslant \left\{
\begin{array}{cc}
  r, &\text{если } r=3,  \\
  r-1, &  \text{если } r>3,
\end{array}
\right.\quad\text{и}
\quad \alpha(x,L)\leqslant s-1.
$$
\end{Lemma}

\begin{proof}\rm Известно, что порядки всех спорадических групп делятся на 2, 3 и 5, и наименьшее простое число, которое может не делить порядок спорадической группы, равно 7. Поэтому $s\geqslant 7.$ Из леммы~\ref{DMPZ} следует, что для любого неединичного элемента $x\in L$ выполнены неравенства
$$
\alpha(x,L)\leqslant 6\leqslant s-1,
$$
и тем самым второе неравенство в утверждении леммы доказано.

 Пусть теперь $|x|$~--- простое число и $r$ делит $|L|$. Если  $|x|=r$, то $\beta_r(x,L)=1$. Поэтому можно считать, что $|x|\ne r$, в частности, $|x|\ne 3$, если $r=3$. Если $\alpha(x,L)\leqslant 3$, то в силу того, что $\beta_r(x,L)\leqslant\alpha(x,L)$, первое неравенство в заключении леммы выполнено. Поэтому мы считаем, что $\alpha(x,L)>3$.
 Мы можем также считать, что $(L,x) \ne (Suz, 3a)$, поскольку в противном случае $3 \leqslant \alpha(x,L)\leqslant 4$ и из того, что $r\ne|x|=3$, получаем $r\geqslant 5$ и $$\beta_r(x,L)\leqslant \alpha(x,L)=4\leqslant r-1.$$

  Таким образом, по лемме~\ref{DMPZ} и с учетом того, что в $M$ имеется два класса сопряженности инволюций $2a$ и $2b$ \cite{atlas}, нам осталось рассмотреть случаи, когда $(L,x)$~--- пара из следующего списка:
$ (J_2, 2a), (Co_2, 2a), (B, 2a), (Fi_{22}, 2a), (Fi_{23}, 2a), (M, 2a), (M, 2b)$. Случаи $(L,x,r)\in\{(Fi_{22}, 2a, 3),(Fi_{23}, 2a, 3)\}$ можно не рассматривать, поскольку в группах Фишера $Fi_{22}$ и $Fi_{23}$ класс сопряженности, обозначенный $2a$, является классом $3$-транспозиций (см. пояснения в~\cite{atlas}), и произведение любых двух некоммутирующих транспозиций из этого класса имеет порядок~$3$. Принимая во внимание неравенство $\beta_r(x,L)\leqslant\alpha(x,L)$, мы можем считать также, что $(L,x,r)$~--- тройка из следующего списка
\begin{multline*}
  (J_2, 2a,3), (Co_2, 2a, 3), (B, 2a,3),  (Fi_{22}, 2a, 5),
  (Fi_{23}, 2a,5), (M, 2a,3), (M, 2b,3).
\end{multline*}

Здесь и далее мы будем пользоваться известным фактом из теории характеров, утверждающим, что для данных элементов $a,b$ и $c$ группы $G$ число $\mathrm{m}(a,b,c)$ пар $(u,v)\in a^G\times b^G$ таких, что $uv=c$, может быть найдено из таблицы характеров по формуле
  $$
  \mathrm{m}(a,b,c)=\frac{|a^G||b^G|}{|G|}\sum\limits_{\chi\in\mathrm{Irr}(G)}\frac{\chi(a)\chi(b)\overline{\chi(c)}}{\chi(1)},
  $$
  см. \cite[упр.~(3.9) на стр.~45]{Isaacs}. В частности, если $G$~--- простая группа,  $|c|=r$ и $
  \mathrm{m}(a,a,c)>0$, то $\beta_r(a,L)\leqslant 2$ (более того, $\beta_r(a,L)\leqslant 2$, если $(|a|,r)=1$). Кроме того, по лемме~\ref{estim} если  $|c|=r$, $\mathrm{m}(a,a,b)>0$ и $\mathrm{m}(b,b,c)>0$,
то $\beta_r(a,L)\leqslant 4$ и т.~д.

Используя таблицы характеров из~\cite{atlas}, рассматриваем\footnote{Вычисления проводились с помощью системы GAP \cite{GAP}.} оставшиеся случаи.
\begin{tabbing}
  $(L,x,r)=(M_{22},2b,5)$: \= $\mathrm{m}(2b,2b,3a)=17060302448280$ \= $\quad\Rightarrow\quad$  \= $\beta_r(x, L)=  2<3=r.$ \kill
  $(L,x,r)=(J_2, 2a)$: \> $\quad\mathrm{m}(2a,2a,3b)=3$ \> $\quad\Rightarrow\quad$ \> $\beta_r(x, L)=  2<3=r.$ \\
  $(L,x,r)=(Co_2, 2a, 3)$: \> $\quad\mathrm{m}(2a,2a,3b)=3$ \> $\quad\Rightarrow\quad$ \> $\beta_r(x, L)=  2<3=r.$ \\
  $(L,x,r)=(B, 2a,3)$: \> $\quad\mathrm{m}(2a,2a,3a)=3$ \> $\quad\Rightarrow\quad$ \> $\beta_r(x, L)=  2<3=r.$ \\
          $(L,x,r)=(Fi_{22}, 2a, 5)$: \> $\quad\mathrm{m}(2a,2a,3a)=3,$ \>                            \> \\
                               \> $\quad\mathrm{m}(3a,3a,5a)=5$  \>    $\quad\Rightarrow\quad$  \> $\beta_r(x, L)\leqslant 4= r-1.$\\
   $(L,x,r)=(Fi_{23}, 2a, 5)$: \> $\quad\mathrm{m}(2a,2a,3a)=3,$ \>                            \> \\
                               \> $\quad\mathrm{m}(3a,3a,5a)=5$  \>    $\quad\Rightarrow\quad$  \> $\beta_r(x, L)\leqslant 4= r-1.$\\
$(L,x,r)=(M, 2a, 3)$: \> $\quad\mathrm{m}(2a,2a,3a)=920808$ \> $\quad\Rightarrow\quad$ \> $\beta_r(x, L)=  2<3=r.$ \\
$(L,x,r)=(M, 2b, 3)$: \> $\quad\mathrm{m}(2b,2b,3a)=17060302448280$ \> $\quad\Rightarrow\quad$ \> $\beta_r(x, L)=  2<3=r.$ \\
\end{tabbing}
     \end{proof}

\begin{Lemma}\label{Outer}
Пусть $L$~--- одна из $26$ спорадических групп, $r\in\pi(L)$~--- нечетное число и $s\notin\pi(L)$. Тогда для любого элемента $x\in \mathrm{Aut}(L)\setminus L$ простого порядка выполнены неравенства $$\beta_r(x,L)\leqslant \left\{
\begin{array}{cc}
  r, &\text{если } r=3,  \\
  r-1, &  \text{если } r>3,
\end{array}
\right.\quad\text{и}
\quad \alpha(x,L)\leqslant s-1.
$$
\end{Lemma}

\begin{proof}\rm
Так как $|\mathrm{Out}(L)|\in\{1,2\}$, элемент $x$~--- инволюция. По известному следствию оригинальной теоремы Бэра~---~Сузуки $x$ инвертирует некоторый элемент $y\in L$ нечетного порядка, и $y$ равен произведению двух элементов, сопряженных с~$x$, откуда $\alpha(x,L)\leqslant2\alpha(y,L)$. Поскольку в диэдральной группе $\langle x,y\rangle$ все инволюции сопряжены и инвертируют любую степень элемента $y$, мы можем считать, что $y$~--- элемент простого порядка.

Как мы заметили при доказательстве леммы~\ref{Inner}, $s\geqslant 7$. По лемме~\ref{DMPZ} либо $\alpha(y,L)\leqslant 3$ и следовательно $$\alpha(x,L)\leqslant2\alpha(y,L)\leqslant 6\leqslant s-1,$$ либо  $(L,y) = (Suz, 3a)$ и $3 \leqslant \alpha(y,L)\leqslant 4$. Но для группы $L=Suz$ наименьшее простое число, не делящее $|L|$, равно $17$. Следовательно, $s\geqslant 17$ и $$\alpha(x,L)\leqslant2\alpha(y,L)= 8<16\leqslant s-1.$$
Второе неравенство в утверждении леммы доказано.

Докажем первое неравенство. Если $\alpha(y,L)\leq 3$ (т.~е. $L\ne Suz$ по лемме~\ref{DMPZ}),  достаточно рассмотреть случаи, когда $r<7$, поскольку в этих случаях $$\beta_r(x,L)\leqslant 2\beta_r(y,L)\leqslant 2\alpha(y,L)\leqslant 6\leqslant r-1.$$ Для $L=Suz$ считаем, что $r\leqslant 7$. Ясно также, что  $|\mathrm{Out}(L)|=2$ в рассматриваемом случае, т.~е.  $$L\in\{M_{12}, M_{22},J_2,J_3,Suz, HS,McL,He, Fi_{22},Fi_{24}',HN, O'N, \}.$$ Поэтому в обозначениях таблиц характеров групп $\mathrm{Aut}(L)$, взятых в \cite{atlas}, считаем, что $(L,x,r)$~--- тройка из следующего списка:
\begin{multline*}
  (M_{12},2c,3),(M_{12},2c,5),(M_{22},2b,3),(M_{22},2b,5),(M_{22},2c,3),(M_{22},2c,5),\\ (J_{2},2c,3), (J_{2},2c,5), (J_{3},2b,3), (J_{3},2b,5), (HS,2c,3),(HS,2c,5), (HS,2d,3),  (HS,2d,5),\\ (McL,2b,3), (McL,2b,5), (He,2c,3),(He,2c,5), (O'N,2b,3), (O'N,2b,5),(HN,2c,3),\\
 (HN,2c,5), (Fi_{22},2d,3),(Fi_{22},2d,5),(Fi_{22},2e,3),(Fi_{22},2e,5), (Fi_{24}',2c,3),(Fi_{24}',2c,5),\\
(Fi_{24}',2d,3),(Fi_{24}',2d,5),
  (Suz,2c,3), (Suz,2c,5), (Suz,2c,7),\\ (Suz,2d,3),(Suz,2d,5), (Suz,2d,7).
\end{multline*}
Тройку $(Fi_{24}',2c,3)$ можно не рассматривать, поскольку $2c$~--- класс $3$-транспозиций в группе Фишера $Fi_{24}$ (см. пояснения в~\cite{atlas}). Разберем каждый из оставшихся случаев, подобно тому, как это сделано в лемме~\ref{Inner}.

\begin{tabbing}
  $(L,x,r)=(M_{22},2b,5)$: \= $\mathrm{m}(3a,3a,5a)=46484685$ \= $\quad\Rightarrow\quad$  \= $\beta_r(x, L)=  2<3=r.$ \kill
    $(L,x,r)=(M_{12},2c,3)$: \> $\quad\mathrm{m}(2c,2c,3b)=18$ \> $\quad\Rightarrow\quad$ \> $\beta_r(x, L)=  2<3=r.$ \\
   $(L,x,r)=(M_{12},2c,5)$: \> $\quad\mathrm{m}(2c,2c,5a)=10$ \> $\quad\Rightarrow\quad$ \> $\beta_r(x, L)=  2<4=r-1.$ \\
    $(L,x,r)=(M_{22},2b,3)$: \> $\quad\mathrm{m}(2b,2b,3a)=3$ \> $\quad\Rightarrow\quad$ \> $\beta_r(x, L)=  2<3=r.$ \\
     $(L,x,r)=(M_{22},2b,5)$: \> $\quad\mathrm{m}(2b,2b,3a)=3,$ \>                            \> \\
                               \> $\quad\mathrm{m}(3a,3a,5a)=500$  \>    $\quad\Rightarrow\quad$  \> $\beta_r(x, L)\leqslant 4= r-1.$\\
      $(L,x,r)=(M_{22},2c,3)$: \> $\quad\mathrm{m}(2c,2c,3a)=9$ \> $\quad\Rightarrow\quad$ \> $\beta_r(x, L)=  2<3=r.$ \\
       $(L,x,r)=(M_{22},2c,5)$: \> $\quad\mathrm{m}(2c,2c,5a)=5$ \> $\quad\Rightarrow\quad$ \> $\beta_r(x, L)=  2<4=r-1.$ \\
        $(L,x,r)=(J_{2},2c,3)$: \> $\quad\mathrm{m}(2c,2c,3b)=18$ \> $\quad\Rightarrow\quad$ \> $\beta_r(x, L)=  2<3=r.$\\
        $(L,x,r)=(J_{2},2c,5)$: \> $\quad\mathrm{m}(2c,2c,3b)=18,$ \>                            \> \\
                               \> $\quad\mathrm{m}(3b,3b,5a)=825$  \>    $\quad\Rightarrow\quad$  \> $\beta_r(x, L)\leqslant 4= r-1.$\\
        $(L,x,r)=(J_{3},2b,3)$: \> $\quad\mathrm{m}(2b,2b,3a)=90$ \> $\quad\Rightarrow\quad$ \> $\beta_r(x, L)=  2<3=r.$\\                         $(L,x,r)=(J_{3},2b,5)$: \> $\quad\mathrm{m}(2b,2b,3a)=90,$ \>                            \> \\
                               \> $\quad\mathrm{m}(3a,3a,5a)=55$  \>    $\quad\Rightarrow\quad$  \> $\beta_r(x, L)\leqslant 4= r-1.$\\
       $(L,x,r)=(HS,2c,3)$: \> $\quad\mathrm{m}(2c,2c,3a)=3$ \> $\quad\Rightarrow\quad$ \> $\beta_r(x, L)=  2<3=r.$\\                           $(L,x,r)=(HS,2c,5)$: \> $\quad\mathrm{m}(2c,2c,3a)=3,$ \>                            \> \\
                               \> $\quad\mathrm{m}(3a,3a,5a)=50$  \>    $\quad\Rightarrow\quad$  \> $\beta_r(x, L)\leqslant 4= r-1.$\\
       $(L,x,r)=(HS,2d,3)$: \> $\quad\mathrm{m}(2d,2d,3a)=75$ \> $\quad\Rightarrow\quad$ \> $\beta_r(x, L)=  2<3=r.$\\                                                   $(L,x,r)=(HS,2d,5)$: \> $\quad\mathrm{m}(2d,2d,5a)=100$ \> $\quad\Rightarrow\quad$ \> $\beta_r(x, L)=  2<4=r-1.$ \\
       $(L,x,r)=(McL,2b,3)$: \> $\quad\mathrm{m}(2b,2b,3a)=810$ \> $\quad\Rightarrow\quad$ \> $\beta_r(x, L)=  2<3=r.$\\                                                   $(L,x,r)=(McL,2b,5)$: \> $\quad\mathrm{m}(2d,2d,5a)=150$ \> $\quad\Rightarrow\quad$ \> $\beta_r(x, L)=  2<4=r-1.$ \\
       $(L,x,r)=(He,2c,3)$: \> $\quad\mathrm{m}(2c,2c,3a)=378$ \> $\quad\Rightarrow\quad$ \> $\beta_r(x, L)=  2<3=r.$\\                                                   $(L,x,r)=(He,2c,5)$: \> $\quad\mathrm{m}(2c,2c,5a)=50$ \> $\quad\Rightarrow\quad$ \> $\beta_r(x, L)=  2<4=r-1.$ \\
       $(L,x,r)=(O'N,2b,3)$: \> $\quad\mathrm{m}(2b,2b,3a)=108$ \> $\quad\Rightarrow\quad$ \> $\beta_r(x, L)=  2<3=r.$\\                                                   $(L,x,r)=(O'N,2b,5)$: \> $\quad\mathrm{m}(2b,2b,5a)=60$ \> $\quad\Rightarrow\quad$ \> $\beta_r(x, L)=  2<4=r-1.$ \\
        $(L,x,r)=(HN,2c,3)$: \> $\quad\mathrm{m}(2c,2c,3a)=3972$ \> $\quad\Rightarrow\quad$ \> $\beta_r(x, L)=  2<3=r.$\\                                                   $(L,x,r)=(HN,2c,5)$: \> $\quad\mathrm{m}(2c,2c,5a)=5250$ \> $\quad\Rightarrow\quad$ \> $\beta_r(x, L)=  2<4=r-1.$ \\
        $(L,x,r)=(Fi_{22},2d,3)$: \> $\quad\mathrm{m}(2d,2d,3c)=3$ \> $\quad\Rightarrow\quad$ \> $\beta_r(x, L)=  2<3=r.$\\                                                    $(L,x,r)=(Fi_{22},2d,5)$: \> $\quad\mathrm{m}(2d,2d,3c)=3,$ \>                            \> \\
                               \> $\quad\mathrm{m}(3c,3c,5a)=32100$  \>    $\quad\Rightarrow\quad$  \> $\beta_r(x, L)\leqslant 4= r-1.$\\
                $(L,x,r)=(Fi_{22},2e,3)$: \> $\quad\mathrm{m}(2e,2e,3a)=13608$ \> $\quad\Rightarrow\quad$ \> $\beta_r(x, L)=  2<3=r.$\\                                                   $(L,x,r)=(Fi_{22},2e,5)$: \> $\quad\mathrm{m}(2e,2e,5a)=80$ \> $\quad\Rightarrow\quad$ \> $\beta_r(x, L)=  2<4=r-1.$ \\
                $(L,x,r)=(Fi_{24}',2c,5)$: \> $\quad\mathrm{m}(2c,2c,3a)=3,$ \>                            \> \\
                               \> $\quad\mathrm{m}(3a,3a,5a)=5$  \>    $\quad\Rightarrow\quad$  \> $\beta_r(x, L)\leqslant 4= r-1.$\\
                $(L,x,r)=(Fi_{24}',2d,3)$: \> $\quad\mathrm{m}(2d,2d,3a)=46484685$ \> $\quad\Rightarrow\quad$ \> $\beta_r(x, L)=  2<3=r.$\\                                                   $(L,x,r)=(Fi_{24}',2d,5)$: \> $\quad\mathrm{m}(2d,2d,5a)=6480$ \> $\quad\Rightarrow\quad$ \> $\beta_r(x, L)=  2<4=r-1.$ \\
       $(L,x,r)=(Suz,2c,3)$: \> $\quad\mathrm{m}(2c,2c,3a)=1620$ \> $\quad\Rightarrow\quad$ \> $\beta_r(x, L)=  2<3=r.$\\                                                   $(L,x,r)=(Suz,2c,5)$: \> $\quad\mathrm{m}(2c,2c,5b)=5$ \> $\quad\Rightarrow\quad$ \> $\beta_r(x, L)=  2<4=r-1.$ \\
       $(L,x,r)=(Suz,2c,7)$: \> $\quad\mathrm{m}(2c,2c,7a)=7$ \> $\quad\Rightarrow\quad$ \> $\beta_r(x, L)=  2<6=r-1.$ \\
       $(L,x,r)=(Suz,2d,3)$: \> $\quad\mathrm{m}(2d,2d,3a)=13608$ \> $\quad\Rightarrow\quad$ \> $\beta_r(x, L)=  2<3=r.$\\                                                   $(L,x,r)=(Suz,2d,5)$: \> $\quad\mathrm{m}(2d,2d,5a)=180$ \> $\quad\Rightarrow\quad$ \> $\beta_r(x, L)=  2<4=r-1.$ \\
       $(L,x,r)=(Suz,2d,7)$: \> $\quad\mathrm{m}(2d,2d,7a)=21$ \> $\quad\Rightarrow\quad$ \> $\beta_r(x, L)=  2<6=r-1.$ \\
  \end{tabbing}
\end{proof}

Теорема~\ref{main} следует из лемм~\ref{Inner} и~\ref{Outer} и неравенства $\beta_s(x,L)\leqslant\alpha(x,L)$. \hfill{$\square$}

 \section{Доказательство теоремы~\ref{Cor}}

 Наряду с теоремой~\ref{main} нам понадобится

 \begin{Lemma} \label{beta_A_n_prop} {\rm \cite[предложение~2]{YRV}}
Пусть  $L=A_n$, $n\geq 5$, $r\leq n$~--- нечетное простое число, и $x\in \mathrm{Aut} (L)$~--- элемент простого порядка. Тогда
\begin{itemize}
\item[$(1)$] $\beta_{r}(x,L)=r-1$, если $x$~--- транспозиция;
\item[$(2)$] справедливо одно из следующих утверждений:
 \begin{itemize}
\item[{\rm (a)}] $\beta_{r}(x,L)\leq r-1$;
\item[{\rm (б)}] $r=3$, $n=6$, $x$~--- инволюция, не лежащая в $S_6$, и   $\beta_{r}(x,L)= 3$.
\end{itemize}\end{itemize}
\end{Lemma}

Доказательство теоремы~{\rm \ref{Cor}}.
Пусть собственное подмножество $\pi$ множества всех простых чисел содержит по крайней мере два различных элемента.  Пусть $r$~--- наименьшее простое число, не принадлежащее~$\pi$, и
$$
m=\left\{
\begin{array}{ll}
  r, & \text{если } r\in\{2,3\}, \\
  r-1, & \text{если } r>3.
\end{array}
\right.
$$
  Допустим, теорема~\ref{Cor} неверна. Тогда $r\geq 3$ по лемме~\ref{2notinpi}, и контрпример $G$ наименьшего порядка к теореме~\ref{Cor} согласно лемме~\ref{red} содержит нормальную подгруппу $L$, изоморфную спорадической или знакопеременной группе, и элемент $x$ простого порядка такие, что $L$ не является $\pi$- или $\pi'$-группой, $G=\langle L,x\rangle$ и любые $m$ сопряженных с $x$ элементов порождают $\pi$-подгруппу в $G$. Кроме того, $C_G(L)=1$, и мы можем, следовательно, считать, что $L\leqslant G\leqslant\mathrm{Aut}(L)$. Пусть $s$~--- наименьший простой делитель $|L|$, не принадлежащий $\pi$. Тогда либо $r$ делит $|L|$ и $s=r$, либо $r$ не делит $|L|$ и $s>r$. По теореме~\ref{main} и лемме~\ref{beta_A_n_prop} выполнено неравенство $$\beta_s(x,L)\leq m,$$ т.~е. существуют элементы ${g_1,\dots,g_m\in L}$ такие, что $|\langle{ x^{g_1},\dots,x^{g_m}}\rangle|$ делится на $s$. Но это противоречит тому, что любые $m$ сопряженных с $x$ элементов порождают $\pi$-под\-группу.

\end{document}